\documentclass[10pt,a4paper]{amsart}

\usepackage{amssymb}
\usepackage{amsthm}
\usepackage{amsfonts}
\usepackage{microtype}
\usepackage{mathrsfs}
\usepackage{graphicx}
\usepackage{color}
\usepackage{latexsym}
\usepackage{rotating}
\usepackage{xspace}
\usepackage[all]{xy}
\usepackage{longtable}
\usepackage{leftidx}
\usepackage{mathtools}
\usepackage{titlesec}
\usepackage{tikz}
\usetikzlibrary{calc}
\usetikzlibrary{arrows}
\usetikzlibrary{decorations.markings}
\usepackage{geometry}

\newtheorem{theorem}{Theorem}[section]
\numberwithin{equation}{theorem}
\newtheorem{lemma}[theorem]{Lemma}

\newtheorem{corollary}[theorem]{Corollary}

\theoremstyle{definition}
\newtheorem{definition}[theorem]{Definition}
\newtheorem{example}[theorem]{Example}

\theoremstyle{conjecture}

\newtheorem{question}[theorem]{Question}
\newtheorem{acknowledgement}{Acknowledgement}

\newcommand{\Ass}{\operatorname{Ass}}
\newcommand{\im}{\operatorname{im}}

\newcommand{\grade}{\operatorname{grade}}

\newcommand{\ara}{\operatorname{ara}}

\newcommand{\cd}{\operatorname{cd}}

\newcommand{\pd}{\operatorname{pd}}

\newcommand{\V}{\operatorname{V}}

\newcommand{\Supp}{\operatorname{Supp}}

\newcommand{\Ann}{\operatorname{Ann}}
\newcommand{\Rad}{\operatorname{Rad}}

\newcommand{\lo}{\longrightarrow}
\newcommand{\fm}{\frak{m}}
\newcommand{\fp}{\frak{p}}

\newcommand{\fa}{\frak{a}}
\newcommand{\fb}{\frak{b}}

\newcommand{\suchthat}{\;\ifnum\currentgrouptype=16 \middle\fi|\;}

\newenvironment{prf}[1][Proof]{\begin{proof}[\bf #1]}{\end{proof}}

\makeatletter
\newcommand{\holim@}[2]{%
  \vtop{\m@th\ialign{##\cr
    \hfil$#1\operator@font holim$\hfil\cr
    \noalign{\nointerlineskip\kern1.5\ex@}#2\cr
    \noalign{\nointerlineskip\kern-\ex@}\cr}}%
}
\newcommand{\holim}{%
  \mathop{\mathpalette\holim@{\rightarrowfill@\textstyle}}\nmlimits@
}
\makeatother

\makeatletter
\def\@secnumfont{\bfseries}
\def\section{\@startsection{section}{1}%
  \z@{.7\linespacing\@plus\linespacing}{.5\linespacing}%
  {\normalfont\Large\bfseries\filcenter}}
\def\subsection{\@startsection{subsection}{2}%
  \z@{.5\linespacing\@plus.7\linespacing}{-.5em}%
  {\normalfont\large\bfseries}}
\makeatother
\begin{document}

\author[K. Divaani-Aazar, A. Ghanbari Doust, M. Tousi and Hossein Zakeri]
{Kamran Divaani-Aazar, Akram Ghanbari Doust, Massoud Tousi\\ and\\ Hossein Zakeri}

\title[Cohomological dimension and relative Cohen-Maculayness]
{Cohomological dimension and relative Cohen-Maculayness}

\address{K. Divaani-Aazar, Department of Mathematics, Alzahra University, Vanak, Post Code 19834, Tehran, Iran-and-School of Mathematics,
Institute for Research in Fundamental Sciences (IPM), P.O. Box 19395-5746, Tehran, Iran.}
\email{kdivaani@ipm.ir}

\address{A. Ghanbari Doust,  Faculty of Mathematical Sciences and Computer, Kharazmi University, Tehran, Iran.}
\email{fahimeghanbary@yahoo.com}

\address{M. Tousi, Department of Mathematics, Shahid Beheshti University, G.C., Tehran, Iran-and-School of Mathematics, Institute for Research
in Fundamental Sciences (IPM), P.O. Box 19395-5746, Tehran, Iran.}
\email{mtousi@ipm.ir}

\address{H. Zakeri, Faculty of Mathematical Sciences and Computer, Kharazmi University, Tehran, Iran.}
\email{hoszakeri@gmail.com}

\subjclass[2010]{13C14; 13C05; 13D45.}

\keywords {Arithmetic rank; cohomological dimension; generalized fractions; local cohomology; relative Cohen-Macaulay module;  system of parameters.}

\begin{abstract} Let $R$ be a commutative Noetherian (not necessary local) ring with identity and $\fa$ be a proper ideal of $R$. We introduce a notion
of $\fa$-relative system of parameters and characterize them by using the notion of cohomological dimension. Also, we present a criterion of relative
Cohen-Macaulay modules via relative system of parameters.
\end{abstract}

\maketitle

\section{Introduction}

Throughout, the word ring stands for commutative Noetherian rings with identity. Consider the following naturally-raised questions:

\begin{question}\label{1.1} Let $\fa$ be a proper ideal of a ring $R$, $M$ a finitely generated $R$-module and $c=\cd(\fa,M)$. Is there a sequence $x_1,
x_2, \ldots, x_c$ of elements in $\fa$ such that $$\cd\left(\fa,M/\langle x_{1}, x_2, \ldots, x_{i} \rangle M\right)=c-i$$ for every $i=1, 2, \ldots, c$?
If yes, how can we characterize such sequences?
\end{question}

Let $R$ be a ring, $\fa$ an ideal of $R$ and $M$ a finitely generated $R$-module with $M\neq \fa M$. Then $M$ is said to be $\fa$-relative Cohen-Macaulay,
$\fa$-RCM, if $\grade(\fa,M)=\cd(\fa,M)$. This notion was introduced by Majid Rahro Zargar and the fourth author in \cite{RZ2} and its study was continued
in \cite{Ra1}, \cite{Ra2}, \cite{Ra3}, \cite{RZ1} and \cite{CH}. Relative Cohen-Macaulay bigraded modules were already introduced and investigated by
Ahad Rahimi; see  \cite{R} and \cite{JR}. Also, the closely related notion of cohomologically complete intersection ideals was examined by Michael
Hellus and Peter Schenzel in \cite{HS}.

\begin{question}\label{1.2} Over a local ring $T,$ a finitely generated $T$-module $N$ is Cohen-Macaulay if and only if every system of parameters of
$N$ is an $N$-regular sequence. Is there an analogue characterization for $\fa$-relative Cohen-Macaulay $R$-modules?
\end{question}

This paper is dealing with the above questions. Although these questions don't look so related in the beginning, surprisingly, they become connected through
a notion of relative system of parameters. Here, we introduce this notion and through investigation of its properties, we answer the above questions.

Let $c:=\cd(\fa,M)$ denote the cohomological dimension of $M$ with respect to $\fa$; i.e. the supermum of the integers $i$ for which $\text{H}_{\fa}^i(M)\neq 0$.
Recall that when $R$ is local with the unique maximal ideal $\fm$ and $\dim_RM=d$, a sequence $x_1, x_2,\ldots, x_d\in \fm$ is called a system of parameters
of $M$ if the $R$-module $M/\langle x_1, x_2,\ldots, x_d\rangle M$ has finite length. This is equivalent to say that $$\Rad\left( \langle x_1, x_2, \ldots,
x_d \rangle +\Ann_RM\right)=\Rad\left(\fm+\Ann_RM \right).$$  We call a sequence $x_{1},x_2, \ldots, x_{c}\in \fa$ an $\fa$-relative system of parameters,
$\fa$-Rs.o.p, of $M$ if $$\Rad\left(\langle x_{1}, x_2, \ldots, x_{c} \rangle+\Ann_{R}M\right)=\Rad\left(\fa+\Ann_{R}M\right).$$

System of parameters appear in many contexts. Especially, Monomial Conjecture on system of parameters of local rings stands for decades until recently
solved by Yves Andr$\acute{e}$; see \cite{An}.
Although over a local ring every finitely generated $R$-module possesses a system of parameters, this is not the case for $\fa$-relative systems of
parameters. It is immediate that $R$ admits an $\fa$-relative system of parameters if and only if $\ara(\fa)=\cd(\fa,R)$. Let $K$ be a field. For a
square-free monomial ideal $\fa$ of a polynomial ring $R=K[x_1,\dots ,x_n]$, it is known that $\cd(\fa,R)=\pd_R\frac{R}{\fa}$; see \cite[Theorem 1]{Ly}.
Characterizing monomial ideals $\fa$ satisfying $\ara(\fa)=\pd_R\frac{R}{\fa}$ has been an active area of research for years; see e.g. \cite{Ba1}, \cite{Ba2}
and \cite{SV}.

Assume that $\fa$ is contained in the Jacobson radical of $R$ and $M$ possesses an $\fa$-Rs.o.p.  We prove that a sequence $\underline{x}=x_1,x_2,\ldots,
x_c\in \fa$ is $\fa$-relative system of parameters of $M$ if and only if $$\cd \left(\fa,M/\langle x_1, x_2, \ldots, x_i \rangle M \right)=c-i$$ for every
$1\leq i\leq c$; see Theorem \ref{2.9}. Also, we show that $M$ is $\fa$-relative Cohen-Macaulay if and only if every $\fa$-relative system of parameters of
$M$ is an $M$-regular sequence if and only if there exists an $\fa$-relative system of parameters of $M$ which is an $M$-regular sequence; see Theorem
\ref{3.2}. These two results yields that if $M$ is $\fa$-RCM and $\underline{x}=x_{1}, x_2, \ldots, x_{c}\in \fa$ is an $\fa$-Rs.o.p of $M$, then
$M/\langle x_{1}, x_2, \ldots, x_{i} \rangle M$ is $\fa$-RCM for every $i=1,\ldots,c$; see Corollary \ref{3.4}.

\section{Question 1.1}

Theorem \ref{2.9} is the main result of this paper. To prove it, we need Lemmas \ref{2.4}, \ref{2.5},  \ref{2.6} and \ref{2.8}. We begin by recalling
some needed definitions.

Let $\fa$ be an ideal of $R$ and $M$ a finitely generated $R$-module. Recall that the arithmetic rank of $\fa$, denoted by $\ara\left(\fa\right)$, is
the least number of elements of $R$ required to generate an ideal with the same radical as $\fa$. Among other things, this paper deals with the local
cohomology modules $$\text{H}_{\fa}^{i}\left(M\right):=\varinjlim \limits_{n\in \mathbb{N}} \text{Ext}_R^i\left(R/\fa^n,M\right); \  i\in \mathbb{N}_0.$$
If $\fb$ is another ideal of $R$ such that the ideals $\fa+\Ann_{R}M$ and $\fb+\Ann_{R}M$ have the same radical, then the
Independence Theorem \cite[Theorem 4.2.1]{BS} yields a natural $R$-isomorphism $\text{H}_{\fa}^{i}\left(M\right)\cong \text{H}_{\fb}^{i}\left(M\right)$
for all $i\in \mathbb{N}_0.$ One easily sees that $\cd\left(\fa,M\right)=-\infty$ if and only if $M=\fa M$. On the other hand, \cite[ Corollary 3.3.3]{BS}
implies that $\cd\left(\fa,M\right)\leq \ara\left(\fa\right)$. In the case $\left(R,\fm\right)$ is a local ring, it is known that $\ara\left(\fm\right)=\dim R=
\cd\left(\fm,R\right)$.

\begin{definition}\label{2.1} Let $M$ be a finitely generated $R$-module and $\fa$ an ideal of $R$ with $M\neq \fa M$.
\begin{enumerate}
\item[i)] Let $c=\cd\left(\fa,M\right)$. A sequence $x_{1},x_2, \ldots, x_{c}\in \fa$ is called $\fa$-{\it relative system of parameters}, $\fa$-Rs.o.p, of $M$ if
$$\Rad\left(\langle x_{1}, x_2, \ldots, x_{c}\rangle+\Ann_{R}M\right)=\Rad\left(\fa+\Ann_{R}M\right).$$
\item[ii)] {\it Arithmetic rank} of $\fa$ with respect to $M$, $\ara\left(\fa,M \right)$, is defined as the infimum of the integers $n\in \mathbb{N}_0$ such that
there exist $x_1, x_2, \ldots, x_n\in R$ satisfying $$\Rad\left(\langle x_{1}, x_2, \ldots, x_{n}\rangle +\Ann_{R}M\right)=\Rad\left(\fa+\Ann_{R}M\right).$$
\end{enumerate}
\end{definition}

Clearly if $x_{1},x_2, \ldots, x_{c}\in R$ is an $\fa$-Rs.o.p of $M$, then for all $t_{1},\ldots,t_{c}\in \mathbb{N}$, every permutation of $x_{1}^{t_{1}},\ldots ,
x_{c}^{t_{c}}$ is also an $\fa$-Rs.o.p of $M$. One may easily check that $\cd\left(\fa,M\right)\leq \ara\left(\fa,M \right)$. Obviously, $\ara\left(\fa,R\right)=
\ara\left(\fa\right)$.

Our first result provides a characterization for existence of relative system of parameters. Although it is an easy observation, we include its proof for the
reader's convenience.

\begin{lemma}\label{2.2} Let $M$ be a finitely generated $R$-module and $\fa$ an ideal of $R$ with $M\neq \fa M$. Then $\fa$ contains an $\fa$-Rs.o.p
of $M$ if and only if $\ara\left(\fa,M \right)=\cd\left(\fa,M\right)$.
\end{lemma}

\begin{prf} Set $c:=\cd\left(\fa,M\right)$. Let $x_{1},x_2, \ldots, x_{c}\in \fa$ be an $\fa$-Rs.o.p of $M$. Then $$\Rad\left(\langle x_{1},x_2, \ldots, x_{c}\rangle+\Ann_RM\right)=\Rad\left(\fa+
\Ann_RM\right),$$ and so $\ara\left(\fa,M \right)\leq c \leq \ara\left(\fa,M \right)$. Thus $\ara\left(\fa,M \right)=c$.

Next, suppose that $\ara\left(\fa,M \right)=\cd\left(\fa,M\right)$. Hence, there are $y_{1}, y_2, \ldots, y_{c}$ in $R$ such that $$\Rad\left(\langle y_{1},\ldots,y_{c}\rangle+\Ann_RM\right)=\Rad\left(\fa+
\Ann_RM\right).$$ There is $n\in \mathbb{N}$ such that $y_{i}^{n}\in \fa+\Ann_RM$ for every $1\leq i\leq c$. So for each $1\leq i\leq c$, there are
$z_i\in \fa$ and $w_i\in \Ann_RM$ such that $y_i^n=z_i+w_i$. Now, $$\Rad\left(\langle z_{1},z_{2}, \ldots, z_{c}\rangle+\Ann_RM\right)=
\Rad\left(\fa+\Ann_RM\right),$$ and so $z_{1},z_{2},\ldots, z_{c}\in \fa$ is an $\fa$-Rs.o.p of $M$.
\end{prf}

This note is also concerned with the special case of the notion of generalized fractions. This notion is described as follows: Let $x_{1}, \ldots, x_{n}$
be a sequence of elements of $R$ and $M$ an $R$-module. Set $$U:=\left \{ \left(x_{1}^{\alpha_{1}}, \ldots, x_{n}^{\alpha_{n}}\right)|\alpha_{1}, \ldots,
\alpha_{n} \in \mathbb{N} \right \}.$$ Then $U$ leads to a module of generalized fractions $U^{-n}M$:\\
For every $r,s\in M$ and $\left(x_{1}^{\alpha_{1}}, \ldots , x_{n}^{\alpha_{n}}\right), \left(x_{1}^{\beta_{1}}, \ldots , x_{n}^{\beta_{n}}\right)\in
U$, we write $$\left(r,\left(x_{1}^{\alpha_{1}}, \ldots , x_{n}^{\alpha_{n}}\right)\right)\sim \left(s,\left(x_{1}^{\beta_{1}},
\ldots , x_{n}^{\beta_{n}}\right)\right)$$ if there exist integers $\delta_i \geq \max\lbrace \alpha_{i}, \beta_{i}\rbrace$; $i=1, \ldots, n$
such that $$x_{1}^{\delta_1-\alpha_{1}}\ldots x_{n}^{\delta_n-\alpha_{n}}r-x_{1}^{\delta_1-\beta_{1}}\ldots x_{n}^{\delta_n-\beta_{n}}s\in
\langle x_{1}^{\delta_1}, \ldots,
x_{n-1}^{\delta_{n-1}}\rangle M.$$ It is easy to verify that $\sim$ is an equivalence relation on  $M\times U$. Then the equivalence class
of an element $\left(r,\left(x_{1}^{\alpha_{1}}, \ldots , x_{n}^{\alpha_{n}}\right)\right)$ is denoted by $\frac{r}{\left(x_{1}^{\alpha_{1}}, \ldots,
x_{n}^{\alpha_{n}}\right)}$ and we let $U^{-n}M$ stand for the set of all equivalence classes of $\sim$. With naturally defined
sum and scalar multiplication, $U^{-n}M$ forms an $R$-module. It is easy to see that $\frac{r}{\left(x_{1}^{\alpha_{1}}, \ldots,
x_{n}^{\alpha_{n}}\right)}\in U^{-n}M$ is zero if and only if there exists an integer $\delta \geq \max\lbrace \alpha_{1},\ldots ,
\alpha_{n}\rbrace$ such that $x_{1}^{\delta-\alpha_{1}}\ldots x_{n}^{\delta-\alpha_{n}}r\in \langle x_{1}^{\delta}, \ldots, x_{n-1}^{\delta}\rangle M$.
For more details see \cite{SZ}.

The next result is very crucial in this paper and may also have applications in other contexts.

\begin{lemma}\label{2.4} Let $\fa=\langle x_{1},\ldots, x_{d}\rangle$ be an ideal of $R$ and $M$ a finitely generated $R$-module.
Then for every $i=1,\ldots, d,$ one has the following exact sequence $$\text{H}_{\fa}^{d-i}\left(\frac{M}{\langle x_{1},\ldots, x_{i-1},
x_{i}\rangle M}\right)\rightarrow \text{H}_{\fa}^{d-i+1}\left(\frac{M}{\langle x_{1}, \ldots, x_{i-1}\rangle M}\right) \stackrel{x_{i}}
\longrightarrow \text{H}_{\fa}^{d-i+1}\left(\frac{M}{\langle x_{1},\ldots, x_{i-1}\rangle M}\right) \rightarrow 0.$$
In particular, there is an exact sequence $$\text{H}_{\fa}^{d-1}\left(M/x_{1}M\right)\longrightarrow \text{H}_{\fa}^{d}\left(M\right) \stackrel{x_{1}}
\longrightarrow \text{H}_{\fa}^{d}\left(M\right) \longrightarrow 0.$$
\end{lemma}

\begin{prf} We first prove the last assertion. Denote $M/x_{1}M$ by $\overline{M}$ and let $-: M\lo \overline{M}$ be the natural
epimorphism. Set $$U:=\left \{ \left(x_{1}^{\alpha _{1}}, x_{2}^{\alpha _{2}}, \ldots, x_{d}^{\alpha _{d}},1\right)| \ \alpha_{1},\ldots,
\alpha_{d}\in \mathbb{N} \right \}$$ and $$V:=\left \{ \left(x_{2}^{\alpha_{2}}, x_{3}^{\alpha_{3}}, \ldots, x_{d}^{\alpha_{d}}, 1\right)|\
\alpha_{2},\ldots,\alpha_{d} \in\mathbb{N} \right \}.$$
Then, by \cite[Remark 2.2]{KSZ}, $\text{H}_{\fa}^{d}\left(M\right)\cong U^{-d-1}M$  and $\text{H}_{\fa}^{d-1}\left(\overline{M}\right)\cong V^{-d}\overline{M}$.
Define $$\varphi:  V^{-d}\overline{M}\lo U^{-d-1}M$$ by $$\varphi\left(\frac{\overline{r}}{\left(x_{2}^{\alpha_{2}}, x_{3}^{\alpha_{3}}, \ldots,
x_{d}^{\alpha_{d}}, 1\right)}\right)
=\frac{r}{\left(x_{1}, x_{2}^{\alpha_{2}}, \ldots, x_{d}^{\alpha_{d}},1\right)}$$ and let $\psi:U^{-d-1}M\lo U^{-d-1}M$ denote the map defined by
multiplication by $x_1$.
It suffices to show that the sequence $$V^{-d}\overline{M}\overset{\varphi}\lo U^{-d-1}M\overset{\psi}\lo U^{-d-1}M\lo 0$$ is exact.
Let $z=\frac{r}{\left(x_{1}^{\alpha_{1}},\ldots,x_{d}^{\alpha_{d}},1\right)} \in U^{-d-1}M  $. Then $z=\frac{x_{1}r}{\left(x_{1}^{\alpha_{1}+1},
x_{2}^{\alpha_{2}}, \ldots, x_{d}^{\alpha_{d}},1\right)}$, and so $\psi$ is surjective.

Clearly, $\im \varphi \subseteq \ker \psi$. Now, let $\psi\left(z\right)=0$. Then, there is an integer $\delta \geq \max\lbrace \alpha_{1},\ldots
,\alpha_{d}\rbrace$ such that $x_{1}^{\delta-\alpha_{1}}\ldots x_{d}^{\delta-\alpha_{d}}x_{1}r \in \langle x_{1}^{\delta},\ldots,x_{d}^{\delta}
\rangle M$.
Hence, $$x_{1}^{\delta+1-\alpha_{1}} x_{2}^{\delta-\alpha_{2}}\ldots x_{d}^{\delta-\alpha_{d}}r=\sum \limits_{i=1}^{d} x_{i}^{\delta}r_{i},$$
where $r_1, \ldots, r_d\in M$. This yields
that
$$\begin{array}{ll}
z&=\frac{r}{\left(x_{1}^{\alpha_{1}},\ldots ,x_{d}^{\alpha_{d}},1\right)}\\ &=\frac{x_{1}^{\delta+1-\alpha_{1}}x_{2}^{\delta-\alpha_{2}}\ldots x_{d}^{\delta-\alpha_{d}}r}{\left(x_{1}^{\delta+1},x_{2}^{\delta},\ldots,x_{d}^{\delta},1\right)}\\
&=\underset{i=1}{\overset{d}\sum }\frac{x_{i}^{\delta}r_{i}}{\left(x_{1}^{\delta +1},x_{2}^{\delta}, \ldots,x_{d}^{\delta},1\right)}\\
&=\frac{x_{1}^{\delta}r_{1}}{\left(x_{1}^{\delta+1},x_{2}^{\delta},\ldots,x_{d}^{\delta},1\right)}\\
&=\frac{r_{1}}{\left(x_{1},x_{2}^{\delta},\ldots,x_{d}^{\delta},1\right)}\\
&=\varphi \left(\frac{\overline{r}_{1}}{\left(x_{2}^{\delta}, x_{3}^{\delta}, \ldots, x_{d}^{\delta}, 1\right)}\right).
\end{array}
$$
So, $\ker \psi \subseteq \im \varphi $.  This completes the proof of the last assertion.

Now, we show the first assertion. Set $N:=M/\langle x_{1},...,x_{i-1}\rangle M$. Then by using the same argument as above, we have the following exact sequence
$$\text{H}_{\fb}^{d-i}\left(N/x_{i}N\right)\longrightarrow \text{H}_{\fb}^{d-i+1}\left(N\right) \stackrel{x_{i}}\longrightarrow \text{H}_{\fb}^{d-i+1}
\left(N\right) \longrightarrow 0,$$ where $\fb:=\langle x_{i},x_{i+1},\ldots, x_{d}\rangle$. This yields our claim.
\end{prf}

\begin{lemma}\label{2.5} Let $M$ be a finitely generated $R$-module and $\fa$ an ideal of $R$ with $M\neq \fa M$. Let $c:=\cd\left(\fa,M\right)$ and
$\underline{x}=x_{1}, x_2, \ldots, x_{c}\in \fa$ be an $\fa$-Rs.o.p of $M$. Then for every $0\leq i\leq c$, the sequence $x_{i+1}, x_{i+2}, \ldots, x_{c}$
forms an $\fa$-Rs.o.p of $M/\langle x_{1},x_2, \ldots, x_{i}\rangle M$, and so $\cd\left(\fa,M/\langle x_{1},x_2, \ldots, x_{i}\rangle M\right)=c-i$.
\end{lemma}

\begin{prf}  We do induction on $i$. The case $i=0$ holds trivially. Next, assume that $i>0$ and the claim holds for $i-1$.  Set $\overline{M}:=M/
\langle x_{1},x_2, \ldots, x_{i-1}\rangle  M$. Then, by the induction hypothesis, $\cd\left(\fa,\overline{M}\right)=c-i+1$. As $$\Rad\left(\fa+\Ann_RM\right)=
\Rad\left(\langle x_{1}, x_2, \ldots, x_{c} \rangle+\Ann_RM\right),$$ Lemma \ref{2.4} yields the exact sequence: $$\text{H}_{\fa}^{c-i}\left(\overline{M}/x_i
\overline{M}\right)\rightarrow \text{H}_{\fa}^{c-i+1}\left(\overline{M}\right) \stackrel{x_{i}}
\longrightarrow \text{H}_{\fa}^{c-i+1}\left(\overline{M}\right) \rightarrow 0.$$
Since $\text{H}_{\fa}^{c-i+1}\left(\overline{M}\right)$ is $\fa$-torsion and $x_i\in \fa$, each element of $\text{H}_{\fa}^{c-i+1}\left(\overline{M}\right)$
is annihilated by some power of $x_i$. Hence, as $\text{H}_{\fa}^{c-i+1}\left(\overline{M}\right)$ is nonzero, the map
$$\text{H}_{\fa}^{c-i+1}\left(\overline{M}\right) \stackrel{x_{i}}\longrightarrow \text{H}_{\fa}^{c-i+1}\left(\overline{M}\right)$$ is not injective. So, $\text{H}_{\fa}^{c-i}\left(\overline{M}/x_i \overline{M}\right)\neq 0$. Consequently, $$c-i\leq \cd\left(\fa,\overline{M}/x_i \overline{M}\right)=
\cd\left(\fa,M/\langle x_{1},x_2, \ldots, x_{i}\rangle M\right).$$

Now, one has the following display of equalities:
$$\begin{array}{ll}
\Rad\left(\fa+\Ann_R\left(M/\langle x_{1},x_2, \ldots, x_{i}\rangle M\right) \right)&=\Rad\left(\fa+\left(\langle x_{1},x_2, \ldots, x_{i}\rangle+\Ann_RM\right)\right)\\
&=\Rad\left(\fa+\Ann_RM\right)\\
&=\Rad\left(\langle x_{1}, x_2, \ldots, x_{c} \rangle+\Ann_RM\right)\\
&=\Rad\left(\langle x_{i+1},  \ldots, x_{c}\rangle+\left(\langle x_{1},x_2, \ldots, x_{i}\rangle+\Ann_RM\right)\right)\\
&=\Rad\left(\langle x_{i+1}, \ldots, x_{c}\rangle+\Ann_R\left(M/\langle x_{1},x_2, \ldots, x_{i}\rangle M\right)\right).
\end{array}
$$
So, $\ara\left(\fa,M/\langle x_{1},x_2, \ldots, x_{i}\rangle M\right)\leq c-i.$ Thus $$\ara\left(\fa,M/\langle x_{1},x_2, \ldots, x_{i}\rangle M\right)=\cd\left(\fa,M/\langle x_{1},x_2, \ldots, x_{i}\rangle M\right)=c-i,$$ and the sequence $x_{i+1}, x_{i+2}, \ldots, x_{c}$ is an $\fa$-Rs.o.p of $M/\langle x_{1},x_2, \ldots, x_{i}\rangle M$.
\end{prf}

Let $\fa$ be an ideal of $R$ and $M, N$ two finitely generated $R$-modules such that $\Supp_RN\subseteq \Supp_RM$. Then, by \cite[Theorem 2.2]{DNT}, $\cd(\fa,N)\leq \cd(\fa,M)$.
In particular if $\Supp_RN=\Supp_RM$, then $\cd(\fa,N)=\cd(\fa,M)$. In the rest of the paper, we shall use this several times without any further comment.

\begin{lemma}\label{2.6} Let $\fa$ be an ideal of $R$ which is contained in the Jacobson radical of $R$
and $x$ an element of $\fa$. Assume that $M$ is a nonzero finitely generated $R$-module with $\ara\left(\fa,M \right)=1$.
If $\text{H}_{\fa}^{1}\left(M/xM\right)=0$, then $\Rad\left(\fa+\Ann_RM\right)=\Rad\left(\langle x \rangle+\Ann_RM\right)$.
\end{lemma}

\begin{prf} Set $\overline{M}:=M/xM$ and assume that $\text{H}_{\fa}^{1}\left(\overline{M}\right)=0$.
As $\fa$ is contained in the Jacobson radical of $R$, it follows that $\overline{M}\neq\fa
\overline{M}$, and so $\cd\left(\fa,\overline{M}\right)\geq 0$. Since $\text{H}_{\fa}^{1}\left(\overline{M}\right)=0$ and $\ara\left(\fa,M \right)=1$,
one deduces that $\cd\left(\fa,\overline{M}\right)=0$. Set $T:=R/\left(\langle x \rangle+\Ann_RM\right)$. Then $\cd\left(\fa,T\right)=\cd\left(\fa,\overline{M}\right)=0$.

There is $y\in R$ such that $$\Rad\left(\fa+\Ann_RM\right)=\Rad\left(\langle y \rangle+\Ann_RM\right).$$  So, $\text{H}_{\fa}^{i}\left(T\right)=\text{H}_{\langle y \rangle}^{i}\left(T\right)$ for every $i\geq 0$. We may and do choose $y$ in $\fa$. By \cite[Remark 2.2.20]{BS}, there is the following exact
sequence $$0 \longrightarrow \text{H}_{\langle y \rangle}^{0}\left(T\right) \longrightarrow T \longrightarrow T_{y} \longrightarrow \text{H}_{\langle y
\rangle}^{1}\left(T\right) \longrightarrow 0,$$ which implies that the natural map $\theta : T\longrightarrow T_{y}$ is surjective. In particular,
there is $t\in T$ such that $\frac{t}{1_T}=\frac{1_T}{y}$, and so $y^{n}\left(yt-1_T\right)=0_{T}$ for some $n\in \mathbb{N}$. As $yt$ belongs to the Jacobson
radical of $T$, $yt-1_T$ is a unite in $T$, and so it follows that $y^{n}\in \langle x \rangle+\Ann_RM$. Thus, $$\Rad\left(\fa+\Ann_RM\right)=\Rad\left(\langle x \rangle+\Ann_RM\right).$$
\end{prf}

A special case of the next result has already been proved by Michael Hellus; see \cite{He2} and \cite[Remark 1.2]{He1}.

\begin{lemma}\label{2.8} Let $\fa$ be a proper ideal of $R$ and $M$ a nonzero finitely generated $R$-module. Let $n\in \mathbb{N}$ be such that
$\cd\left(\fa,M\right)\leq n $ and $x_1,\ldots, x_n \in \fa$. Consider the following conditions:
\begin{enumerate}
\item[i)] $\Rad\left(\langle x_1 \ldots, x_{n} \rangle+\Ann_RM\right)=\Rad\left(\fa+\Ann_RM\right)$.
\item[ii)] The map $\text{H}_{\fa}^{n-i+1}\left(M/\langle x_1 \ldots, x_{i-1} \rangle M\right)\stackrel{x_{i}}\longrightarrow \text{H}_{\fa}^{n-i+1}
\left(M/\langle x_1 \ldots, x_{i-1} \rangle M\right)$ is surjective for all $i=1,\ldots, n$.
\end{enumerate}
Then i) implies ii). Furthermore if $\fa$ is contained in the Jacobson radical of $R$, then i) and ii) are equivalent.
\end{lemma}

\begin{prf} i)$\Rightarrow$ii) It follows by Lemma \ref{2.4}.

ii)$\Rightarrow$i) We do induction on $n$. Assume that $n=1$. Since
$$\Rad\left(\langle x_1 \rangle+\Ann_RM\right)\subseteq \Rad\left(\fa+\Ann_RM\right),$$ it suffices to show that $$\V\left(\langle x_1 \rangle+\Ann_RM\right)\subseteq \V\left(\fa+\Ann_RM\right).$$ Let $\fp\in \V\left(\langle x_1 \rangle+\Ann_RM\right) $.
Then $\fp\in \Supp_RM$. Set $T:=R/\Ann_RM $. Then $$\cd\left(\fa T,T\right)=\cd\left(\fa,M\right)\leq 1.$$
So, $\text{H}_{\fa}^{1}\left(-\right)$ is a right exact endofunctor on the category of $T$-modules and $T$-homomorphisms. By the assumption, the map
$\text{H}_{\fa}^{1}\left(M\right)\stackrel{x_1}\longrightarrow \text{H}_{\fa}^{1}\left(M\right) $ is surjective. Now, we have the following display of $R$-isomorphisms:
$$\begin{array}{ll}
\text{H}_{\fa}^{1}\left(M\right)\otimes_{R} R/\fp &\cong \text{H}_{\fa T}^{1}\left(M\right) \otimes_{T} T/\fp T \\
&\cong \left(\text{H}_{\fa T}^{1}\left(T\right)\otimes_{T} M\right)\otimes_{T} T/\fp T \\
&\cong \text{H}_{\fa T}^{1}\left(T\right)\otimes_{T} M/\fp M \\
&\cong \text{H}_{\fa T}^{1}\left(M/\fp M\right)\\
&\cong \text{H}_{\fa}^{1}\left(M/\fp M\right).
\end{array}$$
This shows that the natural map $$\text{H}_{\fa}^{1}\left(M/\fp M\right)\stackrel{x_1}\longrightarrow \text{H}_{\fa}^{1}\left(M/\fp M\right)$$
is surjective. But $x_1\in \fp$, and so the above map is zero. Thus, $\text{H}_{\fa}^{1}\left(M/\fp M\right)=0$. Since $\fa$ is contained in the Jacobson
radical of $R$ and $\fp\in \Supp_RM$, it turns out that $M/\fp M \neq \fa \left(M/\fp M\right)$, and so $\Gamma_{\fa}\left(M/\fp M\right)\neq 0$.
One has
$$\Rad\left(\Ann_R\left(M/\fp M\right)\right)=\Rad\left(\fp+\Ann_RM\right)=\fp,$$
and so $\Supp_R\left(M/\fp M\right)=\Supp_R\left(R/\fp\right).$ This implies that $$\cd\left(\fa,R/\fp\right)=\cd\left(\fa,M/\fp M\right)=0.$$ Hence
$\Gamma_{\fa}\left(R/\fp\right)\neq 0$, which implies that $\fa \subseteq \fp$, and so $\fp \in \V\left(\fa+\Ann_RM\right)$.

Next, assume that $n>1$ and the case $n-1$ is settled. Set $\overline{M}:=M/x_1M$. As $$\cd\left(\fa,M/\left(0:_{M}x_1\right)\right)\leq \cd\left(\fa,M\right)\leq n ,$$
applying the functor $\text{H}_{\fa}^{n}\left(-\right)$ on the
exact sequence $$0\longrightarrow M/\left(0:_{M}x_1\right) \stackrel{x_1}\longrightarrow M\longrightarrow \overline{M} \longrightarrow 0$$ yields that $\text{H}_{\fa}^{n}\left(\overline{M}\right)$
is a quotient of $\text{H}_{\fa}^{n}\left(M\right) $. So, the map $\text{H}_{\fa}^{n}\left(\overline{M}\right)\stackrel{x_1}\longrightarrow \text{H}_{\fa}^{n}\left(\overline{M}\right)$ is surjective.
But this map is zero, and so $\text{H}_{\fa}^{n}\left(\overline{M}\right)=0 $. Thus, $\cd\left(\fa,\overline{M}\right)\leq n-1$.
Since $$\text{H}_{\fa}^{n-i+1}\left(\frac{M}{\langle x_1 \ldots, x_{i-1} \rangle M}\right)=\text{H}_{\fa}^{\left(n-1\right)-\left(i-1\right)+1}
\left(\frac{\overline{M}}{\langle x_2, \ldots, x_{i-1}\rangle \overline{M}}\right),$$ by the induction hypothesis
$$\Rad\left(\langle x_2, \ldots, x_n \rangle+\Ann_R\left(\overline{M}\right)\right)=\Rad\left(\fa +\Ann_R\left(\overline{M}\right)\right).$$
Now by the argument given in the second paragraph of the proof of Lemma \ref {2.5}, we deduce that
$$\Rad\left(\langle x_1, \ldots, x_n\rangle+\Ann_R\left(M\right)\right)=\Rad\left(\fa+\Ann_R\left(M\right)\right).$$
\end{prf}

Now, we are ready to present the main result of this paper.

\begin{theorem}\label{2.9} Let $\fa$ be an ideal of $R$ which is contained in the Jacobson radical of $R$ and $M$ a nonzero finitely generated $R$-module.
Assume that $c:=\cd\left(\fa,M\right)=\ara\left(\fa,M \right)$ and $x_1,\ldots, x_{c}\in \fa$. Then the following are equivalent:
\begin{enumerate}
\item[i)] $x_1,\ldots, x_{c}$ is an $\fa$-Rs.o.p of $M$.
\item[ii)] The map  $\text{H}_{\fa}^{c-i+1}\left(M/\langle x_1, \ldots, x_{i-1}\rangle M\right) \stackrel{x_{i}}\longrightarrow \text{H}_{\fa}^{c-i+1}
\left(M/\langle x_1,
\ldots, x_{i-1}\rangle M\right)$ is surjective for all $i=1, \ldots, c$.
\item[iii)] $\cd\left(\fa,M/\langle x_{1},x_2, \ldots, x_{i}\rangle M\right)=c-i$ for every $i=1, 2, \ldots, c$.
\end{enumerate}
\end{theorem}

\begin{prf} For $c=0$, there is nothing to prove. So, in the rest of the argument, we assume that $c\geq 1$.

i)$\Leftrightarrow$ii) and i)$\Rightarrow$iii) are immediate by Lemmas \ref{2.8} and \ref{2.5}; respectively.

iii)$\Rightarrow$i) We do induction on $c$. Suppose that $c=1$ and set $\overline{M}:=M/x_1 M$. Then $\ara\left(\fa,M \right)=1$ and $\text{H}_{\fa}^{1}\left(\overline{M}\right)\stackrel{\left(iii\right)}=0$. So, Lemma \ref{2.6}
implies that $$\Rad\left(\fa+\Ann_RM\right)=\Rad\left(\langle x_1 \rangle+\Ann_RM\right).$$ Thus $x_1$ is an $\fa$-Rs.o.p of $M$. Next, suppose
that $c>1$ and the claim holds for $c-1$. One has $\cd\left(\fa,\overline{M}\right)\stackrel{\left(iii\right)}=c-1$ and, for each $2\leq i \leq c$, $$\begin{array}{ll}\cd\left(\fa,\frac{\overline{M}}{\langle x_2, x_3, \ldots, x_{i} \rangle\overline{M}}\right)&=\cd\left(\fa,\frac{M}{\langle x_{1}, x_2,
\ldots, x_{i} \rangle M}\right)\\
&\stackrel{\left(iii\right)}=\cd\left(\fa,M\right)-i\\
&=\cd\left(\fa,M\right)-1-\left(i-1\right)\\
&\stackrel{\left(iii\right)}=\cd\left(\fa,\overline{M}\right)-\left(i-1\right).
\end{array}
$$
Hence by the induction hypothesis, the sequence $x_{2}, x_3, \ldots, x_{c}$ forms an $\fa$-Rs.o.p of $\overline{M}$. Thus
$$\Rad\left(\fa+\Ann_R\overline{M}\right)=\Rad\left(\langle x_2, x_3, \ldots, x_{c} \rangle+\Ann_R\overline{M}\right),$$  which implies that
$$\Rad\left(\fa+\Ann_RM\right)=\Rad\left(\langle x_1, x_{2},  \ldots, x_{c} \rangle+\Ann_RM\right).$$
Therefore, $x_{1},\ldots,x_{c}$ is an $\fa$-Rs.o.p of $M$.
\end{prf}

Next, we record the following immediate conclusion which may be interesting in its own right.

\begin{corollary}\label{2.10}
Let $\left(R,\fm\right)$ be a local ring, $M$ a $d$-dimensional nonzero finitely generated $R$-module and $x_{1},\ldots, x_{d}\in \fm$. Then the following
are equivalent:
\begin{enumerate}
\item[i)] $x_{1},\ldots, x_{d}$ is a system of parameters of $M$.
\item[ii)] The map  $\text{H}_{\fm}^{d-i+1}\left(M/\langle x_{1}, \ldots, x_{i-1}\rangle M\right) \stackrel{x_{i}}\longrightarrow \text{H}_{\fm}^{d-i+1}
\left(M/\langle x_{1},
\ldots, x_{i-1}\rangle M\right)$ is surjective for all $i=1,\ldots, d$.
\end{enumerate}
\end{corollary}

Let $\left(R,\fm\right)$ be a local ring. Next, we will mention two results for system of parameters that their
analogues don't hold for relative system of parameters; see Example \ref{2.11}.

First: Every $R$-regular sequence is a part of a system of parameters of $R$.

Second: Let $M$ be a maximal Cohen-Maculay $R$-module and $A$ be a square matrix of size $n$ with entries in $R$. Let $x_1, x_2, \ldots, x_n$ be a system
of parameters of $M$ and $y_1, y_2, \ldots, y_n\in \fm$ be such that $[y_1, y_2, \ldots, y_n]^T=A[x_1, x_2, \ldots, x_n]^T$. Then by \cite[Theorem]{DR}, $y_1,
y_2, \ldots, y_n$ is a system of parameters of $M$ if and only if the map induced by multiplication by det $A$ from $M/\langle x \rangle M$ to $M/\langle y
\rangle M$ is injective.

\begin{example}\label{2.11}Let $K$ be a field, $R=K[[x,z]]$ and $\fa=\langle x \rangle$. Then $\ara\left(\fa\right)=\cd\left(\fa,R\right)=1$. Set $y:=zx$. Then
\begin{enumerate}
\item[i)] One has $y\in \fa\setminus Z_{R}\left(R\right)$ and $x\notin \Rad\left(\langle y \rangle\right)$.
Hence, $\text{H}_{\fa}^{1}\left(\frac{R}{\langle y \rangle}\right)\neq 0$ by Lemma \ref{2.6}. So, by Lemma \ref{2.5}, $y$ is not an $\fa$-Rs.o.p of $R$.
\item[ii)] The natural map $R/\langle x \rangle\overset{z}\lo R/\langle y \rangle$ is injective, $x$ is an $\fa$-Rs.o.p of $R$
and $[y]=[z][x]$, while $y$ is not an $\fa$-Rs.o.p of $R$.
\end{enumerate}
\end{example}

\section{Question 1.2}

Our main result in this section is Theorem \ref{3.2}. To prove it, we need the following  lemma.

\begin{lemma}\label{3.1} Let $\fa=\langle x_{1},\ldots,x_{n}\rangle$ be an ideal of $R$ and $M$ a finitely generated $R$-module with $\fa M\neq M$.
Set $g:=\grade\left(\fa,M\right)$. Then
\begin{enumerate}
\item[i)] $\fa$ can be generated by elements $y_{1},\ldots,y_{n}$ such that $y_{i_{1}}, \ldots, y_{i_{h}}$
forms an $M$-regular sequence for all $i_{1},\ldots,i_{h}$ with $1\leq i_{1}<\cdots<i_{h}\leq n, h\leq g$.
\item[ii)] If $\fa$ is contained in the Jacobson radical of $R$ and $g=n$, then $x_{1},\ldots,x_{n}$ forms an $M$-regular sequence.
\end{enumerate}
\end{lemma}

\begin{prf} i) Follows by \cite[Theorem 125 (b)]{Ka}.

ii)  Follows by \cite[Theorem 129]{Ka}.
\end{prf}

Here is the right place to bring the following immediate corollary of Lemma \ref{2.2}.

\begin{corollary}\label{2.3} Let $M$ be a finitely generated $R$-module and $\fa$ an ideal of $R$ with $M\neq \fa M$. Then the following are equivalent:
\begin{enumerate}
\item[i)] $M$ is $\fa$-RCM and it possesses an $\fa$-Rs.o.p.
\item[ii)] $\grade\left(\fa,M\right)=\ara\left(\fa,M \right)$.
\end{enumerate}
\end{corollary}

\begin{prf} It is clear by  Lemma \ref {2.2} and the inequality $\grade\left(\fa,M\right)\leq \cd\left(\fa,M\right)\leq \ara\left(\fa,M \right)$.
\end{prf}

\begin{theorem}\label{3.2} Let $M$ be a finitely generated $R$-module and $\fa$ an ideal of $R$ with $\ara\left(\fa,M \right)=\cd\left(\fa,M\right)$.
Consider the following conditions:
\begin{enumerate}
\item[i)] $M$ is $\fa$-RCM.
\item[ii)] Every $\fa$-Rs.o.p of $M$ is an $M$-regular sequence.
\item[iii)] There exists an $\fa$-Rs.o.p of $M$ which is an $M$-regular sequence.
\end{enumerate}
Then i) and iii) are equivalent. Furthermore if $\fa$ is contained in the Jacobson radical of $R$, all three conditions are equivalent.
\end{theorem}

\begin{prf} Set $c:=\cd\left(\fa,M\right)$.

i)$\Rightarrow$iii) Let  $y_{1}, y_2, \ldots, y_{c}\in \fa$ be an $\fa$-Rs.o.p of $M$ and set $J:=\langle y_{1},y_2, \ldots, y_{c}\rangle$. Then $$\Rad\left(J+\Ann_RM\right)=\Rad\left(\fa+\Ann_RM\right),$$ and so $$\grade\left(J,M\right)=\grade\left(\fa,M\right)=c.$$
By Lemma \ref{3.1} i), there exist $x_{1}, x_2, \ldots, x_{c}\in R$ which forms an $M$-regular sequence and $J=\langle x_{1}, x_2, \ldots, x_{c}
\rangle$. Now, $x_{1}, x_2, \ldots, x_{c}$ is our desired $\fa$-Rs.o.p of
$M$.

iii)$\Rightarrow$i) Let $\underline{z}=z_{1},\ldots,z_{c}\in \fa$ be an $\fa$-Rs.o.p of $M$ which is an $M$-regular sequence. Then $$\ara\left(\fa,M
\right)\leq c\leq \grade\left(\fa,M\right)\leq \ara\left(\fa,M \right).$$ So, $M$ is $\fa$-RCM by Corollary \ref{2.3}.

ii)$\Rightarrow$iii) It is obvious.

i)$\Rightarrow$ii) Let $\underline{z}=z_{1}, z_2, \ldots, z_{c}\in \fa$ be an $\fa$-Rs.o.p of $R$ and set $J:=\langle z_{1},z_2, \ldots, z_{c}\rangle$.
Then
$$\Rad\left(J+\Ann_RM \right)=\Rad\left(\fa+\Ann_RM\right),$$ and so
$$\begin{array}{ll}
\grade \left(J,M \right)&=\grade \left(J+\Ann_RM,M \right)\\
&=\grade \left(\fa + \Ann_RM , M \right)\\
&=\grade \left(\fa,M \right)\\
&=c.
\end{array}$$
Thus Lemma \ref{3.1} ii) yields that $\underline{z}$ is an $M$-regular sequence.
\end{prf}

The following example shows that in Theorem \ref{3.2}, the assumption that $\fa$ is contained in the Jacobson radical of $R$ is necessary.

\begin{example}\label{3.3}
Let $K$ be a field. Consider the ring $R=K[x,y,z]$ and let $\fa=\langle x,y,z \rangle$. It is clear that $R$ is $\fa$-RCM. We can see that $\fa=\langle
y\left(1-x\right), z\left(1-x\right), x\rangle$, so that $y\left(1-x\right), z\left(1-x\right), x$ is an $\fa$-Rs.o.p of $R$. But $y\left(1-x\right),
z\left(1-x\right), x$ is not an $R$-regular sequence.
\end{example}

Next, we record the following corollary of Theorem \ref{3.2}.

\begin{corollary}\label{3.4}  Let $\fa$ be an ideal of $R$ which is contained in the Jacobson radical of
$R$. Let $M$ be an $\fa$-RCM $R$-module and $\underline{x}=x_{1}, x_2, \ldots, x_{c}\in \fa$ an $\fa$-Rs.o.p
of $M$. Then $M/\langle x_{1}, x_2, \ldots, x_{i} \rangle M$ is $\fa$-RCM for every $i=1,\ldots,c$.
\end{corollary}

\begin{prf} Set $\overline{M}:=M/\langle x_{1}, x_2, \ldots, x_{i} \rangle M$. By Theorem \ref{3.2}, $x_{1}, x_2, \ldots, x_{c}$ is an $M$-regular sequence,
and so $x_{i+1}, x_{i+2}, \ldots, x_{c}$ is an $\overline{M}$-regular sequence.  On the other hand, by Lemma \ref{2.5},  the sequence $x_{i+1},x_{i+2}, \ldots,
x_{c}$ is an $\fa$-Rs.o.p of $\overline{M}$.  Applying Theorem \ref{3.2} again implies that $\overline{M}$ is $\fa$-RCM.
\end{prf}

\begin{acknowledgement}
The authors thank Sara Saeedi Madani for introducing them to the references \cite{Ba1}, \cite{Ba2} and \cite{SV}.
\end{acknowledgement}


\end{document}